\documentclass[11pt,notitlepage]{amsart} 
\usepackage{amssymb,array,delarray,latexsym} 

\newlength{\sh}
\newlength{\jmr}
\newlength{\jfc}
\newlength{\bernd}
\newlength{\jones}
\newlength{\mati}
\newlength{\sil}
\newlength{\koi}
\newtheorem{jst}{The JST Theorem}

\newtheorem{cor}{Corollary}

\newtheorem{dfn}{Definition}

\newtheorem{main}{Main Theorem} 

\newtheorem{rem}{Remark}	
\newtheorem{ex}{Example}

\newcommand{\blah}{\hspace{1in}} 
\newcommand{\pspa}{{$\mathbf{PSPACE}$}} 
\newcommand{\am}{{$\mathbf{AM}$}} 

\newcommand{\nc}{{$\mathbf{NC}$}}

\newcommand{\eps}{\varepsilon}

\newcommand{\cE}{\mathcal{E}}

\newcommand{\oE}{{\bar{E}}}

\newcommand{\cO}{\mathcal{O}}
\newcommand{\supp}{\mathrm{Supp}}
\newcommand{\conv}{\mathrm{Conv}}
\newcommand{\thth}{{\underline{\mathrm{th}}}}

\newcommand{\st}{{\underline{\mathrm{st}}}}

\newcommand{\Q}{\mathbb{Q}}
\newcommand{\R}{\mathbb{R}}
\newcommand{\C}{\mathbb{C}}
\newcommand{\N}{\mathbb{N}}
\newcommand{\Z}{\mathbb{Z}}

\newcommand{\area}{\mathrm{Area}}

\newcommand{\res}{\mathrm{Res}}

\newcommand{\Zn}{\Z^n}
\newcommand{\Zm}{\Z^m}
 
\newcommand{\Qn}{\Q^n} 

\newcommand{\Rm}{\R^m}

\newcommand{\Cn}{\C^n}

\renewcommand{\qed}{$\blacksquare$}
\newcommand{\cM}{{\mathcal{M}}}
\newcommand{\cH}{\mathcal{H}}
\newcommand{\aM}{\cM^{\mathrm{ave}}}

\newcommand{\cT}{\mathcal{T}}
\newcommand{\bO}{\mathbf{O}}
\newcommand{\vol}{\mathrm{Vol}}

\newenvironment{mymatrix}{\begin{array}{*{10}{c}}}{\end{array}}

\begin{document}

\title{On the Complexity of  
Diophantine Geometry in Low Dimensions 
(Extended Abstract)} 

\author{J. Maurice Rojas}
\thanks{ Submitted for publication. 
This research was partially funded by a Hong Kong CERG Grant. } 
\address{Department of Mathematics\\
City University of Hong Kong\\
83 Tat Chee Avenue\\
Kowloon, HONG KONG\\ 
{\tt mamrojas@math.cityu.edu.hk}\\ 
{\tt http://www.cityu.edu.hk/ma/staff/rojas} } 


\date{October 7, 1998} 

\begin{abstract} 
We consider the average-case complexity of some otherwise  
undecidable or open Diophantine problems. More precisely, 
we show that the following two problems can be solved in \pspa: 
\begin{enumerate} 
\renewcommand{\theenumi}{\Roman{enumi}} 
\item{Given polynomials $f_1,\ldots,f_m\!\in\!\Z[x_1,\ldots,x_n]$ defining a 
variety of dimension $\leq\!0$ in $\Cn$, find all solutions 
in $\Zn$ of $f_1\!=\cdots=\!f_m=0$. } 
\item{For a given polynomial $f\!\in\!\Z[v,x,y]$ defining an irreducible 
nonsingular non-ruled surface in $\C^3$, decide the sentence 
$\exists v \; \forall x \; \exists y \; f(v,x,y)\!\stackrel{?}{=}\!0$, 
quantified over $\N$. } 
\end{enumerate} 
Better still, we show that the truth of 
the {\bf Generalized Riemann Hypothesis (GRH)} implies that  
detecting roots in $\Qn$ for the polynomial systems in 
problem (I) can be done via a two-round Arthur-Merlin protocol, i.e., well 
within the second level of the polynomial hierarchy. (Problem (I) is, of 
course, undecidable without the dimension 
assumption.) The decidability of problem (II) was previously unknown. 
Along the way, we also prove new complexity and size bounds for solving 
polynomial systems over $\C$ and $\Z/p\Z$. A 
practical point of interest is that the aforementioned 
Diophantine problems should perhaps be avoided in the construction 
of crypto-systems.  
\end{abstract} 

\maketitle 

\section{Introduction and Main Results}
\label{sec:intro} 
The negative solution of Hilbert's Tenth Problem 
\cite{undec,hilbert10} has all but dashed earlier hopes of solving 
large polynomial systems over the integers. However, an 
immediate positive consequence is the creation of a 
rich and diverse garden of hard problems with potential applications in 
complexity theory, cryptology, and logic. 
Even more compelling is the question of where the boundary to 
decidability lies. 

From high school algebra we know that detecting roots in $\Q$ (or $\Z$ 
or $\N$) for polynomials in $\Z[x_1]$ is tractable. 
However, in \cite{jones9}, Jones showed that detecting roots 
in $\N^9$ for polynomials in $\Z[x_1,\ldots,x_9]$ is already undecidable. 
Put another way, this means that detecting a positive integral point on a 
general algebraic hypersurface of (complex) dimension $8$ is 
undecidable. 

It then comes as quite a shock that decades of number theory 
still haven't settled the complexity of the analogous question for 
algebraic varieties of dimension $1$ through $7$. In fact, even the case 
of plane curves remains a mystery:\footnote{In particular, the major 
``solved'' special cases so far have only extremely ineffective 
complexity and height bounds. (See, 
e.g., the introduction and references of \cite{big}.)} 
As of late 1998, the decidability of 
detecting a root in $\N^2$, $\Z^2$, or even $\Q^2$, for an 
{\bf arbitrary} polynomial in $\Z[x_1,x_2]$, is still completely open.
 
\subsection{Dimension Zero} 
We will thus go one dimension lower\footnote{We use the natural convention that $\dim Z\!:=\!-1$ when $Z=\emptyset$.} in order to prove a useful 
result.  
\begin{main}
\label{main:riemann}
Suppose $F\!:=\!(f_1,\ldots,f_m)$ is a system of polynomials 
in $\Z[x_1,\ldots,x_n]$ and let $Z$ be the zero set of $F$ 
in $\Cn$. Assume further that $\dim Z\!\leq\!0$. Then 
\begin{enumerate}
\item{We can find all the roots of $F$ in $\Zn$ within \pspa. } 
\item{The truth of the Generalized Riemann Hypothesis implies 
that deciding $Z\!\cap\!\Q^n\!\stackrel{?}{=}\!\emptyset$ is in \am.}
\end{enumerate}  
\end{main}
\begin{rem}
Recall that 
$\mathbf{NP}\!\cup\!\mathbf{BPP}\!\subseteq\!\mathbf{AM}\!\subseteq\!
\mathbf{RP^{NP}}\!\cap\!\mathbf{coNP^{NP}}\!\subseteq\!\mathbf{\Pi}_2
\!\subseteq\!\cdots\!\subseteq\!\mathbf{PH}\!\subseteq\!\mathbf{PSPACE}
\!\subseteq\!\mathbf{EXPTIME}$ \cite{lab,papa}. Also, it will later follow 
easily (cf.\ remark \ref{rem:final}) that the rational analogue of assertion 
(1) can be done within $\mathbf{EXPTIME}$ (or polynomial time for 
fixed $n$). 
\end{rem} 
While one can derive assertion (1) more or less directly (from, say, 
\cite{renegar} or \cite{bpr}), the explicit 
sequential and parallel complexity bounds we give   
in section \ref{sec:proof1} (cf.\ remark \ref{rem:final}) are the best to 
date and do {\bf not} follow from earlier work. Also, assertion (2) presents a 
new arithmetic analogue of a recent result of Koiran \cite{hnam} stating 
that the truth of GRH implies that detecting a {\bf complex} root can 
also be done within \am, with no restriction on $\dim Z$. 
\begin{rem} 
Deciding $\dim Z\!\stackrel{?}{\leq}\!0$ can be done in \pspa \,  
via another algorithm of Koiran \cite{koiran}. In fact, the 
truth of GRH implies that this decision problem can be done within 
\am \, as well \cite{koiran}.  
\end{rem} 
\begin{rem} 
If one fixes the monomial term structure of $F$ and picks the coefficients 
randomly, then it follows easily from the theory of resultants 
\cite{gkz94,introres} that $m\!\geq\!n \Longrightarrow F$ will have only 
finitely many roots in $\Cn$ on average. This can be made completely explicit 
via, say, J.\ Schwartz' well-known result \cite{schwartz} on randomized 
verification of polynomial identities. 
\end{rem}  
\begin{rem}
We will use the classical Turing machine \cite{papa} as our computational 
model, along with the bit-wise (resp.\ sparse) encoding for 
rational numbers (resp.\ polynomials in $\Z[x_1,\ldots,x_n]$) 
\cite{bss}, throughout. So, for example, the size of an integer 
$c$ will be $1+\lceil\log_2(|c|+1)\rceil$. Also, all 
quantification symbols will be 
understood to range over the positive integers $\N$. 
\end{rem} 

An interesting corollary of Main Theorem \ref{main:riemann} is the 
following: 
\begin{cor} 
\label{cor:uni} 
Following the notation of Main Theorem \ref{main:riemann}, assume 
$n$ is {\bf fixed}. Then we can find all the roots of $F$ in 
$\Zn$ within $\mathbf{NC}_3$. \qed 
\end{cor}
For $n\!=\!1$, this complements two older results: (a) the famous 
result of Lenstra, Lenstra, and Lovasz \cite{lll} that all the rational 
roots of $f$ can be found in polynomial sequential time, and (b) 
Neff's result \cite{neff} that the roots of $f$ in $\C$ can be 
approximated (to some input precision) in $\mathbf{NC_3}$. 
We have thus obtained higher-dimensional Diophantine analogues of these 
results. 

The proof of Main Theorem \ref{main:riemann} is based on the 
following number-theoretic result relating root-finding over the 
fields $\Q$ and $\Z/p\Z$. 
\begin{main} 
\label{main:start} 
Following the notation of Main Theorem \ref{main:riemann}, let 
$r$ be the number of roots of $F$ in $\Qn$ counting multiplicities. 
Also define $N_F(x)$ to be 
the {\bf total} number of roots of $F$ in $\Z/p\Z$, counting 
multiplicities, as $p$ runs 
through all primes $\leq\!x$. Finally, let $\pi(x)$ denote 
the number of primes $\leq\!x$. Then the truth of GRH implies that  
\[\left|\frac{N_F(x)}{\pi(x)}-r\right|\leq 
\frac{A}{\sqrt{x}}[C_F\log x + \cM(E)(\log x)^2],\] 
where $A$ is an effective and absolute constant, and $C_F$ is a 
quantity polynomial in $S(\oE)$ and the size of $F$. 
\end{main} 
This result may be of independent interest to number 
theorists, as well as complexity theorists and numerical analysts. 
We also note that while Main Theorem \ref{main:start} deals with 
using reduction mod $p$ to count roots over $\Q$, other results, such as 
\cite[Thm.\ 8]{hnam} and \cite[Thm.\ 4.1]{peter}, use reduction mod 
$p$ to determine the existence of roots in $\C$. Our techniques can 
be used to improve these latter results as well, and this 
will be pursued in a forthcoming paper.  
\begin{rem} 
The main results we describe in the next two subsections 
are of a more technical nature, so the reader more interested 
in the conceptually simple \am \, algorithm above can skip directly 
to the beginning of section \ref{sec:proof1}. 
\end{rem}  

\subsection{Dimensions One and Two} 
To reconsider the complexity of detecting 
integral points on varieties of dimension $\geq\!1$, one can consider 
more subtle combinations of quantifiers to facilitate finding  
decisive results. For example, Matiyasevich and Julia Robinson have 
shown \cite{matrob,jones81} that sentences 
of the form $\exists u \; \exists v \; \forall x \; \exists y \; 
f(u,v,x,y)\!\stackrel{?}{=}\!0$ are already undecidable. 

However, the decidability of sentences 
of the form $\exists v \; \forall x \; \exists y \; 
f(v,x,y)\!\stackrel{?}{=}\!0$ 
was an open question until recently: in \cite{big} it was shown that 
these sentences can be decided by a Turing machine, 
once the input $f$ is suitably restricted. Roughly speaking, 
deciding the prefix $\exists\forall\exists$ is equivalent to determining 
whether an algebraic surface has a slice (parallel to the 
$(x,y)$-plane) densely peppered with positive integral 
points.  The ``exceptional'' $f$ not covered by the algorithm of 
\cite{big} form a very slim subset of $\Z[v,x,y]$. 

We will further improve this result by showing that 
under similarly mild input restrictions, $\exists\forall\exists$ 
can in fact be decided in singly exponential time and 
parallelized effectively. 
To make this more precise, let us write any $f\!\in\!\Z[v,x,y]$ 
as $f(v,x,y)\!=\!\sum c_av^{a_1}x^{a_2}y^{a_3}$, where 
the sum is over certain $a\!:=\!(a_1,a_2,a_3)\!\in\!\Z^3$. 
We then define the {\bf support} of $f$ as 
$\mathbf{Supp(f)}:=\{a \; | \; c_a\!\neq\!0 \}$. The {\bf Newton 
polytope} of $f$, $\mathbf{Newt(f)}$, is then just the convex 
hull\footnote{That is, the smallest convex set in $\R^3$ containing 
$\supp(f)$.}  $\conv(\supp(f))$. When we say that a statement involving a set of parameters $\{c_1,\ldots,c_N\}$ is true {\bf generically}, we will mean that 
the statement holds for all $(c_1,\ldots,c_N)\!\in\!\C^N$ outside of some {\bf a priori fixed} algebraic hypersurface. 
\begin{main}
\label{main:pepper} 
Fix the Newton polytope $P$ of a polynomial $f\!\in\!\Z[v,x,y]$  
and let $\rho$ denote the projection mapping $\R^3$ onto the 
$(\hat{e}_2,\hat{e}_3)$-plane. Suppose further that $\rho(P)$ has at 
least one lattice point in its interior. Then, for a generic choice of 
coefficients depending only on $P$, we can decide 
$\exists v \; \forall x \; \exists y \; f(v,x,y)\!\stackrel{?}{=}\!0$ 
within \pspa.
\end{main} 
The generic choice above is clarified further in section 
\ref{sec:proof2}. In particular, we will see that the ``average'' 
polynomial in $\Z[v,x,y]$ with moderately large support defines an 
irreducible, nonsingular, non-ruled surface in $\C^3$.  

It is interesting to note that the exceptional case to our above 
algorithm judiciously contains an extremely hard number-theoretic 
problem: determining the existence of a point in $\N^2$ on an 
algebraic plane curve. (Indeed, we will see later that 
$\Z[v,y]$ lies in our exceptional locus.) We also point out that 
the problem of computing the size of the {\bf largest} positive 
integral point on an algebraic plane curve is closely related to 
determining whether an algebraic surface possesses {\bf any} integral 
point: It was recently shown in \cite{big} that under certain assumptions, 
the decidability of the latter problem implies the uncomputability 
of the former function. 

Aside from a geometric trick, the proof of Main Theorem 2 relies 
on essentially the same tools as the proof of Main Theorem 1:  
Both proofs make use of the toric resultant \cite{gkz94,introres,gcp} 
and a recent perturbation trick specifically tailored for degenerate  
sparse polynomial systems \cite{gcp}. New complexity and 
size estimates on polynomial system solving 
over $\C$ form a crucial key step. We now describe these new 
bounds.  

\subsection{New Bounds for Solving Over $\C$} 
\label{sub:bounds} 

In order to rigourously state our results for polynomial system 
solving over $\C$, 
we will introduce some necessary notation: First note that the 
notions of support and Newton polytope extends naturally to 
polynomials in $\Z[x_1,\ldots,x_n]$. We then say that $F$ is an 
$\mathbf{m\times n}$ {\bf polynomial system with support contained 
in} $\mathbf{E}$, whenever $F$ is as stated in Main Theorem 1, 
$E\!=\!(E_1,\ldots,E_m)$, and $\supp(f_i)\!\subseteq\!E_i$ for 
all $i$.  

An important geometric invariant for $m\times m$ systems of equations 
is $\cM(E)$ --- the {\bf mixed volume} 
\cite{buza,gk94,isawres,ewald,mvcomplex} of the convex hulls 
of $E_1,\ldots,E_m$. A trick we will use to solve general 
$m\times n$ systems (when $m\!\geq\!n$) is to include additional 
points in the $E_i$ so that $\cM(E)\!>\!0$. 

For $(m+1)\times m$ systems we will instead let 
$E\!:=\!(E_1,\ldots,E_m)$ and $\oE\!:=\!(E_1,\ldots,E_{m+1})$. 
We also point out the following two important 
complexity-theoretic parameters: 
\begin{dfn} 
\label{dfn:first} 
For any $(m+1)$-tuple $\oE$ of finite subsets of $\Zm$, define 
$\mathbf{R(\oE)}\!:=\!\sum^{m+1}_{i=1}\cM(E_1,\ldots,
\widehat{E_i},\ldots,E_{m+1})$ and let  
$\mathbf{S(\oE)}\!=\!\cO(\sqrt{m}e^m\aM_{\oE})$, 
where $\aM_{\oE}$ is the average value of $\cM(\cE)$ as $\cE$ ranges over 
all $m$-tuples $(\cE_1,\ldots,\cE_m)$ with $\cE_i\!\in\{E_1,\ldots,E_{m+1}\}$ for all $i$. 
\end{dfn}
While the precise definition of $S(\oE)$ depends on the 
efficiency of a particular class of algorithms described in 
\cite{gcp}, the preceding asymptotic bound will suffice for our purposes. 
\begin{rem}
\label{rem:dense} 
An important extreme class of $F$ is the {\bf dense} case: This 
occurs when, for all $i$, all monomial terms up to some fixed degree occur 
in $f_i$. In this case, the best known complexity bounds for 
existential quantifier elimination over $\C$ are polynomial in 
$D_\Pi$ and \begin{tiny}$\left(\begin{mymatrix} D_\Sigma+1\\ n 
\end{mymatrix}\right)$\end{tiny}, where 
$D_\Pi$ and $D_\Sigma$ are respectively the product and sum of the 
total degrees of the $f_i$ \cite{ckl89,ierardi,fitch,gcp}.  
\end{rem} 
\begin{rem}
\label{rem:upper}  
Our complexity and size bounds will instead be polynomial in the 
invariants $\cM(E)$, $R(\oE)$, and $S(\oE)$. In particular, we 
point out that $\cM(E)\!\leq\!R(\oE)\!\leq\!S(\oE)$. 
Furthermore, $\cM(E)\!\leq\!D_\pi$, $R(\oE)\!\leq\!(n+1)D_\Pi$, and 
$S(\oE)\!\!\leq$\begin{tiny}$\left(\begin{mymatrix} D_\Sigma+1\\ n 
\end{mymatrix}\right)$\end{tiny}$<\!(\frac{eD_\Sigma}{n})^n$, with 
equality if the 
Newton polytopes are the same as those of the dense case. Even better, 
for sparse polynomial systems, our invariants are usually dramatically smaller than thes last three upper limits (cf.\ example \ref{ex:dense}). 
\end{rem}  

To bound the size of the roots of very general polynomial 
systems, we will need one final invariant. 
\begin{dfn}
Following the notation of definition \ref{dfn:first}, 
let $c,k\!\in\!\N$ and define 
\[ H(c,k,\oE)\!:=\!2^{\cM(E)}\left(2m^2\cM(E)^2\right)^{m\cM(E)}
\left(c\sqrt{k}\right)^{S(\oE)-\cM(E)}.\]   
\end{dfn} 

Let $\Delta\!:=\!\conv(\{\bO,\hat{e}_1,\ldots,\hat{e}_m\})$, 
where $\bO\!\in\!\Rm$ denotes the origin and $\hat{e}_i\!\in\!\Rm$ is the 
$i^\thth$ standard basis vector. In what follows, $\cO^*(T)$ means 
$\cO(T\log^r T)$ for some constant $r\!>\!0$. 
\begin{main} 
\label{main:solve}
Following the notation and hypotheses of Main Theorem 1, suppose 
further that $m\!\geq\!n$, every $f_i$ has at most $k$ monomial terms, 
and that all the coefficients of $F$ have absolute value $\leq\!c$. Also, for 
all $i\!\in\!\{1,\ldots,m\}$ define $E_i$ to be the union of 
$\{\bO,\hat{e}_i\}$ and the support of $f_i$.   
Finally, let $E_{m+1}\!\:=\!\Delta\cap\Zm$. Then we 
can find univariate polynomials $h,h_1,\ldots,h_n\!\in\!\Z[t]$ with the 
following properties: 
\begin{enumerate}
\addtocounter{enumi}{-1}
\item{The degrees of $h,h_1,\ldots,h_n$ are bounded 
above by $\cM(E)$.} 
\item{Given any root $\theta\!\in\!\C$ of $h$, define 
$\gamma(\theta)\!:=\!(h_1(\theta),\ldots,h_n(\theta))$. Then 
the set of points $\{\gamma(\theta)\}_{h(\theta)=0}$ contains all 
the roots of $F$ in $\Cn$.} 
\item{The coefficients of $h$ have absolute value 
bounded above by $H(c,k,\oE)$ } 
\end{enumerate}
In particular, the size of the coefficients of $h$ 
is polynomial in $S(\oE)$ and $L_F$, where $L_F$ denotes the size of 
$F$. Furthermore, $h,h_1,\ldots,h_n$ can be found deterministically within 
{\bf sequential} time 
$\cO^*(nm^3\cM(E)^3S(\oE)^{2.376}\min\{R(\oE)^2,S(\oE)\}L_F)$, 
or {\bf parallel} time $\cO(\log \{S(\oE)L_F\})$ using 
\mbox{$\cO^*(nm^3\cM(E)^3S(\oE)^{2.376}\min\{R(\oE)^2,S(\oE)\})$ proc-} 
essors. 
In particular, $h,h_1,\ldots,h_n$ can be found in \nc \, when 
$n$ is fixed. 
\end{main} 
\begin{rem} 
Choosing extra points to add to $\supp(F)$ (and the choice of Newton polytope 
for $f_{m+1}$) can be optimized for a given monomial term structure without 
too much difficulty (see \cite[remark 2]{gcp}).  
\end{rem} 
The complexity bounds above are the best deterministic bounds to date for 
{\bf general}\footnote{The 
algorithm above can easily be generalized to the case $\dim Z\!\geq\!1$ 
without changing the complexity bounds,   
via the techniques of \cite{gcp}. } 
univariate reduction over $\C$. In particular, 
our bounds work for over-determined systems of equations and 
are much more sensitive to the sparse encoding of the input 
than earlier bounds (cf.\ example \ref{ex:dense}). For instance, 
our bounds above considerably improve those of \cite{ierardi,fitch}, which 
were stated solely in terms of the degrees of the input 
polynomials.\footnote{We point out, however, that these papers 
deal with a more general problem than Main Theorem \ref{main:solve}. }  
Furthermore, the exponents above can be significantly lowered if 
randomization is allowed \cite{gcp}. 

As an almost immediate corollary, we obtain new and extremely 
general bounds on the size 
of the roots of $F$: 
\begin{main} 
\label{main:height} 
Following the notation and hypotheses of Main Theorems 
\ref{main:riemann} and \ref{main:solve}, let $(x_1,\ldots,x_n)\!\in\!\C$ 
be any root of $F$,  fix any $j\!\in\!\{1,\ldots,n\}$, and now let $E_{m+1}$ 
be the line segment from the origin $\bO$ to $\hat{e}_j$. Then either 
$x_j\!=\!0$ or $H(c,k,\oE)^{-1}\!\leq\!|x_j|\!\leq\!H(c,k,\oE)$. 
In particular, the size of the roots of $F$ is polynomial in $S(\oE)$ 
and the size of $F$.  
\end{main}
\noindent
The above gap theorem generalizes earlier bounds due 
to Malajovich-Mu\~noz \cite{gregthesis} and Canny \cite{cannyphd}. 
In particular, while our bound is a bit coarser in the dense   
case, it is significantly better for certain sparse systems (cf.\ 
example \ref{ex:dense}). Furthermore, earlier bounds assumed 
$m\!=\!n$ and made various other nondegeneracy assumptions ---  
such as no multiple roots --- even in the case $m\!=\!n\!=\!1$.
 
Main Theorems \ref{main:riemann}, \ref{main:pepper}, \ref{main:solve}, 
and \ref{main:height} are respectively proved in sections 
\ref{sec:proof1}, \ref{sec:proof2}, \ref{sec:proof3}, and \ref{sec:proof4}. 
Section \ref{sec:proof4} also contains the proof of Main Theorem  
\ref{main:start}.   
However, we will first give some examples before starting our main proofs. 

\section{Algorithmic Examples} 
\label{sec:back} 
We begin with a brief illustration and synopsis of the 
algorithm from Main Theorem \ref{main:solve}. We then conclude with an 
example of the benefits of the sensitivity to sparsity (or ``monomial 
sensitivity'') in our complexity and size bounds. 

Further details on the univariate reduction algorithm below appear in section 
\ref{sec:proof3}, and the case $m\!=\!n$ is described at length in \cite{gcp}.  
\begin{ex} {\bf ($\mathbf{m\!=\!n\!=\!2}$)} 
\label{ex:2by2} \\
Consider the bivariate polynomial system $F=(1+2x-2x^2y-5xy+x^2+3x^3y,$\\ 
$2+6x-6x^2y-11xy+4x^2+5x^3y)$.
Letting $E$ be the support of $F$, the reader can easily 
verify\footnote{For $n=2$, there is the simple 
formula $\cM(E)=$\\ $\area(\conv(E_1+E_2))-\area(\conv(E_1))
-\area(\conv(E_2))$. Also, both polynomials are divisible by $x+1$.} 
that $\bO\!\in\!E_1\cap E_2$,  
$\cM(E)\!=\!4$, and that the only roots of $F$ are the points 
$\{(1,1),(\frac{1}{7},\frac{7}{4})\}$ and the line 
$\{-1\}\!\times\!\C$.  So this example illustrates a slightly 
more general situation than our $\dim Z\!=\!0$ restriction. 

We now highlight the main points of our algorithm: First, via combinatorial 
means \cite{isawres}, we construct a {\bf toric resultant matrix}, $M_\oE$. 
The entries of this sparse and highly-structured matrix (of size 
$S(\oE)\times S(\oE)$) will be 
the coefficients of $F$, as well as a few extra parameters. 

The coefficients of the polynomials $h,h_1,h_2$ are found  
by solving a collection of linear systems of equations. 
The coefficients for these systems are in turn derived from the computation 
of determinants of various specializations of $M_\oE$.  
The number of linear systems and determinants, as well as their sizes, 
accounts for the complexity bound of Main Theorem \ref{main:solve}. 
Note in particular that $R(\oE)\!=\!4+4+4\!=\!12$ and that $S(\oE)$ 
can be taken to be $17$ in the case at hand. (The corresponding matrix 
$M_\oE$ appears explicitly in \cite{gcp}.)  

After an application of {\tt Maple}, one at last obtains that  
\[ h(t)  =  -153+120t+1540t^2+1600t^3+448t^4,\]  
$h_1(t) = -\!$\small $\!\frac{11762}{7511}$ \normalsize $+\!$ 
                \small $\!\frac{19150}{22533}\!$ \normalsize $\!t \ +\!$ 
                \small $\!\frac{114736}{22533}\!$ \normalsize $\!t^2 \ +\!$ 
                \small $\!\frac{7264}{3219}\!$ \normalsize $\!t^3$, and  
$h_2(t) = -\!$\small $\!\frac{5881}{7511}$ \normalsize $+\!$ 
           \small $\!\frac{32108}{22533}\!$ \normalsize $\!t \ +\!$
           \small $\!\frac{57368}{22533}\!$ \normalsize $\!t^2 \ +\!$
           \small $\!\frac{3632}{3219}\!$ \normalsize $\!t^3$.     

Since $h(t)$ factors as 
\[(2t+3)(28t+51)(2t+1)(4t-1),\] 
we immediately obtain 
a set of points lying in $Z$ (including all the isolated roots) 
by substituting $\{-\frac{3}{2},-\frac{51}{28},-\frac{1}{2}, 
\frac{1}{4}\}$ into the pair 
$(h_1(t),h_2(t))$. 
Note also that the roots of $h$ are exactly $-\frac{1}{2}\zeta_1-\zeta_2$, 
as $(\zeta_1,\zeta_2)$ ranges over a finite set of roots of $F$, at 
least  one in each irreducible component of $Z$. 
\end{ex} 

Our next example shows how older complexity and size bounds 
(which were stated solely in terms of degrees of polynomials) 
may be too pessimistic for certain sparse systems. 
\begin{ex} {\bf (Well Directed Spikes)} 
\label{ex:dense}
Consider the system of equations $F$ defined by
\[ a_{1,1}+a_{1,2}x_1+\cdots+a_{1,n}x_{n-1} \blah \]
\[ +c_{1,1}(x_1\cdots x_n)+ \cdots + c_{1,d}(x_1\cdots x_n)^d =0 \]
\[ \vdots \] 
\[a_{n,1}+a_{n,2}x_1+\cdots+a_{n,n}x_{n-1} \blah \] 
\[ +c_{n,1}(x_1\cdots x_n)+ \cdots + c_{n,d}(x_1\cdots x_n)^d =0.\] 
In this case, the Newton polytopes are all equal to a single  
``spike'' and we can actually pick the Newton polytope of $f_{n+1}$ 
to be $\conv(\{\bO,\hat{e}_1,\ldots,\hat{e}_{n-1},\sum \hat{e}_j\})$. So 
via the basic properties of the mixed volume, it is a 
routine exercise to check that $\cM(E)\!=\!d$ and $R(\oE)\!=\!(n+1)d$. 
Furthermore, note that our example is equivalent (via the 
change of coordinates $y\!=\!(x_1,\ldots,x_{n-1},x_1\cdots x_n)$) 
to a system of total degree $1$ in $y_1,\ldots,y_{n-1}$ and of degree 
$d$ in $y_n$. It then follows from \cite{gkz94,gcp} that we 
can build a sufficiently compact resultant matrix so that 
$S(\oE)\!=\!R(\oE)$. Main Theorem \ref{main:solve} thus tells 
us that we can reduce solving $F$ to a univariate problem within 
deterministic sequential time $\cO^*(n^{7.376}d^{6.376})$. 

Let us now see how previous methods fare on our sparse 
example. Remark \ref{rem:dense} tells us that older algorithms took sequential 
time polynomial in $D_\Pi$ and a rather large binomial coefficient. 
These two parameters respectively reduce to $n^nd^n$ and 
$\Omega(n^n(\frac{d}{e})^n)$ (by Stirling's estimate). So our 
methods clearly exploit sparsity to much greater advantage. 

As for root size bounds, the best previous bound for dense 
systems of {\bf identical} degrees \cite{cannyphd} specializes 
to $(3cdn)^{d^nn^{n+1}}$. Roughly speaking, 
this means that this older bound needs $\cO(bd^nn^{n+1})$ bits to specify 
the size of the roots of $F$, assuming the coefficients of $F$ used 
$b$ bits to start with. On the other hand, our new bound from 
Main Theorem \ref{main:height} tells us that we need only 
$\cO(nd(b+\log(nd))$ bits. So it appears that sparsity 
is also quite helpful for root size bounds. 

Generating infinite families of such examples is 
easy, simply by picking Newton polytopes which are $n$-dimensional, 
but ``long'' in a suitable fixed direction.  
\end{ex} 

\section{The Proof of Main Theorem \ref{main:riemann}}
\label{sec:proof1}
We will first present a proof of assertion (2), since 
the corresponding algorithm is much simpler than our 
more intricate $\mathbf{PSPACE}$ algorithm. 

\noindent 
{\bf Proof of Assertion (2):} 
Consider the following algorithm: 
\begin{itemize} 
\item[{\bf Step 1}]{ Pick positive integers $s$ and $t$ as follows: 
Let $t$ be just large enough so that
\[ \frac{A}{\sqrt{t}}[C_F\log t+\cM(E)(\log t)^2]\!<\!\frac{1}{4}\]  
and $t\!>\!4n\cM(E)$. Then define $s$ to be just small enough so that 
$\frac{t}{s}\!>\!\frac{\pi(t)}{4}$. } 
\item[{\bf Step 2}]{ Pick a (uniformly) random integer 
in $j\!\in\!\{1,\ldots,\lceil\frac{t}{s}\rceil-1\}$ and 
define $x_j$ and $x_{j+1}$ to be $js$ and $(j+1)s$ respectively. 
(If $j+1\!=\!\lceil\frac{t}{s}\rceil$ then define $x_{j+1}\!:=\!t+1$.) }  
\item[{\bf Step 3}]{ Via an $\mathbf{NP}$ oracle, find 
if there is a prime $p\!\in\![x_j,x_{j+1})$ such that 
$F$ has a solution in $\Z/p\Z$. } 
\item[{\bf Step 4}]{ If such a prime exists, conclude that 
$F$ has a root in $\Qn$. Otherwise, declare that there 
is no root in $\Qn$. }   
\end{itemize} 
This algorithm clearly at least resembles a two-round Arthur-Merlin protocol, 
so let us confirm its correctness. 

First note if $F$ has a root in $\Qn$, then $F$ has a root in $\Z/p\Z$ for 
all but finitely many primes $p$. (The exceptional primes must divide the 
denominator of some coordinate of a rational root of $F$, so there can be 
at most $n\cM(E)$ such primes.) 
So in this case, by our choice of $t$, $F$ will have a root in $\Z/p\Z$ 
for more than $\frac{3}{4}$ of the primes $p\!\in\![1,t+1)$. So, 
by our choice of $s$, our algorithm has a probability greater than 
$\frac{3}{4}$ of succeeding in this case. The role of the size of $s$ is 
detailed further in the next case: 

Suppose now that $F$ has no roots in $\Qn$. Call a prime for 
which $F$ has a root in $\Z/p\Z$ a {\bf bad} prime. Then by Main Theorem 
\ref{main:start} and our choice of $t$, strictly less than 
$\frac{1}{4}$ of the primes $p\!\in\![1,t+1)$ are bad. Also, 
by the box-principle, the probability that the interval $[x_j,x_{j+1})$ 
contains a bad prime is less than $\frac{1}{4}$ (thanks to our choice of 
$s$). So the probability of failure in this case is strictly less than 
$\frac{1}{4}$. 

To conclude, recall that $\log S(\oE)$ is polynomial in the 
size of $F$. So the number of bits necessary to express any prime above 
is polynomial in the size of $F$. Also note that assuming GRH, $\pi(x)$ 
(which is asymptotic to $\frac{x}{\log x}$ \cite{weinberger}) can 
be approximated within a factor of $1+\eps$ in time polynomial in 
$\frac{1}{\eps}$ and the size of $x$ \cite{lago}. (Furthermore, 
we can certainly compute an integer at least as large as $S(\oE)$ in 
polynomial time via remark \ref{rem:upper}.) Thus, the time 
needed to compute $\{x_j,x_{j+1}\}$, the number of necessary 
random bits, and the number of bits of any integer in $[x_j,x_{j+1})$, are 
all polynomial in the size of $F$. 

Finally, via \cite{miller}, the truth of GRH implies that an integer 
$p\!\in\![x_j,x_{j+1})$ can be verified to be a prime in polynomial time. 
So via \cite{cohen}, a putative root of $F$ mod $p$ can indeed be verified 
in polynomial time. So we indeed need only one call to an $\mathbf{NP}$ 
oracle. Our algorithm is thus indeed an \am \, algorithm. \qed 

\noindent 
{\bf Proof of Assertion (1):}  
First note that if $m\!<\!n$ then $Z\!\neq\!\emptyset \Longrightarrow 
\dim Z\!\geq\!1$, by the well-known facts on intersections of 
complex hypersurfaces \cite{mumford}. So we can safely assume that 
$m\!\geq\!n$. 
 
Following the notation of the proof of Main Theorem \ref{main:solve}, 
now pick $f_{m+1}\!=\!u_0+u_ix_i$ and run the algorithm from 
Main Theorem \ref{main:solve} (cf.\ remark \ref{rem:special}) to solve for 
the $i^\thth$ coordinates of all the roots of $F$ within accuracy 
$\frac{1}{3}$. By Main Theorem \ref{main:solve}, and combining 
with Neff's algorithm for \nc \ univariate solving over $\C$ 
\cite{neff}, this computation can be done for all $i\!\in\!\{1,\ldots,n\}$ 
within \pspa. (Note that we are implicitly using the fact that 
polynomial parallel time with singly exponentially work is in 
$\mathbf{PSPACE}$ \cite[pg.\ 398]{papa}.) 

By taking combinations of the coordinates thus found, 
we obtain a set $R\!\subset\!\C^n$ of at most $\cM(E)^n$ putative approximate 
solutions of $F$. In particular, by the choice of accuracy we've made, 
the integral roots of $F$ are contained in the set 
$R'\!:=\!\{ ([x_1],\ldots,[x_n]) \; | \; (x_1,\ldots,x_n)\!\in\!T \}$, 
where $[x]$ denotes the nearest integer to $x$. 
The integral roots of $F$ then certainly lie in $R'$. They can 
then be identified within \pspa, using the most naive parallel algorithm 
for polynomial evaluation, devoting at worst polynomially many 
processors to each element of $R'$. The quantity 
$\cM(E)$ is at worst singly exponential in the coordinates 
of the $E_i$ (cf.\ remark \ref{rem:upper}), so we are done. \qed  

\begin{rem}
\label{rem:final} 
A more precise complexity bound for algorithm above is the 
following: \[\cO^*(n^2m^3\cM(E)^3S(\oE)^{2.376}\min\{R(\oE)^2,S(\oE)\}L_F)\]
sequential time (via fast approximate univariate factorization over 
$\C$ \cite{binipan}) or $\cO^*(\log^3 \{S(\oE)L_F\})$ parallel time 
using $\cO^*(nm^3\cM(E)^3S(\oE)^{2.376}\min\{R(\oE)^2,S(\oE)\})$ processors. 
(So we obtain $\mathbf{NC}_3$ for fixed $n$.) 
Furthermore, it easy to see that we can substitute the 
factoring algorithm of \cite{lll}, in place of Neff's algorithm, to find all 
the roots of $F$ in $\Qn$ with $\mathbf{EXPTIME}$. For fixed 
$n$ this clearly reduces to polynomial time. 
\end{rem} 

\begin{rem} 
Presumably, Main Theorem \ref{main:riemann} continues to hold under the 
weaker condition that the {\bf real} dimension of $Z$ 
is at most zero. One route of proof is an extension of Main Theorem 
\ref{main:height} to the {\bf real} isolated roots of $F$, 
and this will be pursued in later work. 
\end{rem} 

\section{Genus Zero Varieties and the Proof of Main Theorem 
\ref{main:pepper}} 
\label{sec:proof2} 
In what follows, we will make use of some basic algebraic geometry.  
A more precise description of the tools we use can be found in 
\cite{big}. Also, we will always use 
{\bf geometric} (as opposed to arithmetic) genus for  
algebraic varieties \cite{hart}. 

Let us begin by clarifying the genericity condition of Main Theorem \ref{main:pepper}. Let $Z$ be the zero set of $f$. 
What we will actually require of $f$ (in addition to the 
assumptions on its Newton polytope) is that $Z$ be an 
irreducible nonsingular surface of positive genus. That $Z$ is 
irreducible and nonsingular for a generic choice of coefficients 
follows from Bertini's theorem \cite{mumford}. That $Z$ also has 
positive genus generically follows from a result of Khovanskii 
\cite{kho78}. (His result actually implies that for generic coefficients 
and generic $v_0$, $f(v_0,x,y,)\!=\!0$ defines a curve of positive 
genus. It is then impossible for $Z$ to generically have 
genus zero.) Since the intersection of any two open Zariski-dense 
sets is open and dense, we indeed have that our hypothesis 
occurs generically. 

Now note that from the classification of algebraic surfaces 
\cite{beau}, $Z$ has positive genus $\Longrightarrow Z$ is non-ruled. 
In particular, this means that there can only be finitely 
many $v_0$ such that the ``slice'' $Z\cap\{v\!=\!v_0\}$ contains a 
curve of genus zero. (By the nonsingularity of $Z$ and Bertini's theorem 
again.) Note that by Siegel's Theorem 
\cite{wow}, $\forall x\;  \exists y\;  f(v_0,x,y)\!=\!0 
\Longrightarrow Z\cap\{v\!=\!v_0\}$ contains a curve of 
genus zero. So assuming one can decide the prefix $\forall\exists$, 
finding genus zero slices for the large family of $f$ above gives 
us a way to decide the prefix $\exists\forall\exists$. 

The preceding assumption is true, in spades, thanks to the following result: 
\begin{jst}
\cite{jones81,schinzel,tungcomplex}
The quantifier prefix $\forall\exists$ is
decidable in sequential time polynomial in $\deg f$ and singly 
exponential in the size of $f$. More explicitly, 
given $f\!\in\!\Z[x,y]$, we have that
$\forall x\; \exists y\; f(x,y)\!=\!0$
iff all of the following conditions hold:
\begin{enumerate}
\item{The polynomial $f$ factors into the form
$f_0(x,y)\prod^k_{i=1}(y-f_i(x))$ where $f_0(x,y)\!\in\!\Q[x,y]$
has {\bf no} zeroes in the ring $\Q[x]$, and for all $i$,
$f_i\!\in\!\Q[x]$ and the leading coefficient of $f_i$ is positive.}
\item{$\forall x\!\in\!\{1,\ldots,x_0\} \; \exists
y\!\in\!\N$ such that $f(x,y)\!=\!0$, where $x_0\!=\!\max\{s_1,\ldots,s_k\}$,
and for all $i$, $s_i$ is the sum of the squares of the coefficients
of $f_i$.}
\item{Let $\alpha$ be the least positive integer such that
$\alpha f_1,\ldots,\alpha f_k\!\in\!\Z[x]$ and set 
$g_i\!:=\!\alpha f_i$ for all $i$. 
Then the {\bf union} of the solutions of the following $k$ congruences
\begin{eqnarray*}
g_1(x) & \!\equiv & 0 \ \mod \ \alpha \\
 & \vdots & \\
g_k(x) & \!\equiv & 0 \ \mod \ \alpha 
\end{eqnarray*}
is {\bf all} of $\Z/\alpha\Z$. \qed }
\end{enumerate}
\end{jst}
\begin{rem}
The JST Theorem can be strengthened slightly in the
following way: one can replace $\alpha$ in condition (3)
with {\bf any} positive integer $\alpha'$ such that 
$\alpha'f_1,\ldots,\alpha'f_k\!\in\!\Z[x]$. Also, 
the techniques of \cite{cohen} can easily support 
the stated complexity bound. 
\end{rem} 

\noindent 
{\bf Proof of Main Theorem \ref{main:pepper}: } It follows 
easily from the Hurwitz genus formula for curves \cite{sil} 
that $\forall x\;  \exists y\;  f(v_0,x,y)\!=\!0
\Longrightarrow Z\cap\{v\!=\!v_0\}$ defines a curve with a singular 
component. The set of such $v_0\!\in\!\C$ is of course finite (by Bertini's 
theorem again), since $Z$ was assumed to be nonsingular. 

So our algorithm is the following: Find those $v_0\!\in\!\N$ 
for which $Z\cap\{v\!=\!v_0\}$ is singular, and then  
solve the corresponding instances of $\forall\exists$. 
If any instance is true, then our original sentence was 
true. Otherwise, our original sentence was false. 

Finding this set of $v_0$ is easily done within polynomial 
sequential time using the Jacobian criterion for singularity 
\cite{mumford}: simply find those positive integers $v_0$ for which 
the system of equations $(f(v_0,x,y),\frac{\partial f(v_0,x,y)}{\partial x},\frac{\partial f(v_0,x,y)}{\partial 
y})$ has a solution $(x,y)\!\in\!\C^2$. This can be done 
easily via Main Theorem \ref{main:solve}. 
In particular, the number of eligible $v_0$ is polynomial in the degree 
of $f$. (Polynomial in $\vol(P)$ in fact). Furthermore, we can simply solve 
to within accuracy $\frac{1}{3}$ to isolate the $v_0\!\in\!\N$ (if any).  

To conclude, we simply note that the only part of 
the JST theorem which can not be implemented in polynomial sequential 
time, via the results we've introduced and quoted so far, is part (2). For 
this part we then simply use (singly) exponentially many processors (one for 
each $x\!\in\!\{1,\ldots,x_0\}$) to finish in constant additional time. 
Since $\mathbf{P}$ (and constant parallel time with singly 
exponential work) is contained in $\mathbf{PSPACE}$ 
\cite[pg.\ 398]{papa}, we are done. \qed 

\begin{rem}
Although a result weaker than Main Theorem \ref{main:solve} would 
have sufficed, an immediate corollary of our proof 
is that the parallel time sufficient to decide 
$\exists\forall\exists$ is near-heptic in the 
volume of the Newton polytope $P$.  
\end{rem}
\begin{rem} 
Note that if $f\!\in\!\Z[v,y]$ then $Z$ is a ruled 
surface in $\C^3$. From another point of view, the hypothesis of 
Main Theorem 
\ref{main:pepper} is violated since $\rho(P)$ is contained in a line 
segment. Deciding $\exists\forall\exists$ for this case then 
reduces to deciding $\exists\exists$, which we've already observed 
is very hard. Nevertheless, Alan Baker has conjectured that this 
problem is decidable \cite[section 5]{jones81}. 
\end{rem} 
\begin{rem} 
The complexity of deciding whether a given surface is ruled 
is an open problem. (Although one can check certain instances  
in \pspa, as described above.) It is also interesting to note 
that finding explicit parametrizations of {\bf rational} surfaces (a 
special class of ruled surfaces) appears to be decidable. Evidence 
is provided by an algorithm of Josef Schicho which, while still 
lacking a termination proof, seems to work well in practice \cite{schicho}. 
\end{rem} 

\section{Determinants Galore and the Proof of Main Theorem \ref{main:solve}}
\label{sec:proof3} 
Let us first concentrate on an important special case. 

\noindent 
{\bf The Case $\mathbf{m\!=\!n}$:} 
This special case of our result, save for the parallel time and coefficient 
bounds, reduces to Main Theorem 1 of \cite{gcp}. Note in particular 
that the sequential time bound for this case simplifies 
to $\cO^*(n^4\cM(E)^3S(\oE)^{2.376}\min\{R(\oE)^2,S(\oE)\}L_F)$. 
There is also the issue of converting arithmetic cost to bit cost, 
so this is how the factor of $L_F$ (not present in \cite{gcp}) 
appears here. We will soon see that the structure of our algorithm 
permits a relatively simple conversion between these two different 
costs via, say, the techniques of \cite{binipan}.  

We will now prove the parallel complexity bound, illustrating 
our algorithm along the way. What we will describe is essentially 
an outline of an algorithm detailed further in \cite{gcp}: 

\begin{itemize} 
\item[{\bf Step 0}]{Compute the toric resultant matrix 
$M_\oE$ corresponding to the support $\oE$. }  
\item[{\bf Step 1}]{Upon suitably specializing $u_1,\ldots,u_m$ 
to generic integers, and picking a generic polynomial 
system $F^*$ with support contained in $E$, define\footnote{In reality, 
$\cH$ will be a divisor of a variant of this determinant. This accounts 
for the copious amount of interpolation below. See \cite{gcp} 
for further details.} $\cH(s,u_0)$ to be 
$\det M_\oE$. }
\item[{\bf Step 2}]{Let $h(t)$ be the coefficient of the 
term of $\cH(s,t)$ of lowest degree in $s$. }  
\item[{\bf Step 3}]{Compute $h_1,\ldots,h_n$ by a slightly 
more complicated combination of interpolation and  
specialized determinants.} 
\end{itemize} 
(The toric resultant $\res_\star(\cdot)$ appears as follows: $\det M_\oE$ 
is actually a multiple of $\res_\oE(\cdot)$ evaluated at the polynomial system 
$(F,f_{m+1})$, where $f_{m+1}\!=\!u_0+u_1x_1+\cdots+u_mx_m$.)  

Step (0) is a preprocessing step with complexity dominated by the  
remainder of our algorithm. The proof follows easily from 
the theory of \cite{gkz94,isawres} (via the Cayley trick, 
the well-known algorithms for simplicial subdivisions \cite{preparata}, 
and the mixed-subdivision algorithm for computing $M_\oE$),  
but would be too great a 
digression to include in this abstract. So we consider our 
remaining steps.

The determination of a generic $F^*$ in Step (1) can be done in 
work well-dominated by the remainder of our algorithm, by 
Main Theorem 3 of \cite{gcp}. (The parallelization of 
finding $F^*$ is quite straightforward.) 

As for the rest of Steps (1) and (2), a lucky choice of 
$u_1,\ldots,u_n$, would reduce 
our work to evaluating $\cO(n\cM(E)\min\{R(\oE)^2,S(\oE)\})$ 
determinants of size $\cO^*(S(\oE))\times 
\cO^*(S(\oE))$ and one evaluation/interpolation problem of size 
$\cO^*(S(\oE))$. In the absence of such luck, one can instead 
be deterministic and evalute $\cO^*(n^3\cM(E)^3\min\{R(\oE)^2,S(\oE)\})$ 
such determimants and solve $\cO^*(\cM(E)^2)$ such interpolation 
problems. (Via a technique of \cite[section 5.2]{gcp}, this 
results in a choice of $u_i$ with absolute values at 
most $\cO(2^nn^{2n}\cM(E)^{2n})$.) 
In particular, it easily follows from 
recent parallel matrix algorithms, e.g., \cite{binipan}, that 
Steps (1) and (2) take overall (arithmetic) parallel time $\cO(\log S(\oE))$ 
using $\cO^*(n^3\cM(E)^3S(\oE)^{2.376}\min\{R(\oE)^2,S(\oE)\})$ processors. 
In terms of the Turing (bit) model, \cite{binipan} (and the 
structure of our algorithm) tells us that 
we need only include an additional factor of $L_F$. So the 
overall parallel time is $\cO(\log \{S(\oE)L_F\})$ with the same number 
of processors.  

By \cite{gcp} it then follows that Step (3)  
amounts to $n$ repetitions of what we just did for Step (2), 
and this work can be done in parallel. We thus at last arrive 
at a parallel complexity bound of $\cO(\log \{S(\oE)L_F\})$ time using 
$\cO^*(n^4\cM(E)^3S(\oE)^{2.376}\min\{R(\oE)^2,S(\oE)\})$ processors. 

We now analyze the size of the coefficients of $h$. Here, 
we will consider a fundamental special case. The proof of 
the general case follows easily from the special case via 
known bounds on the coefficients of factors of multivariate polynomials, 
e.g., \cite{mignotte}.

So let us make the following assumption: the determinant of 
$M_\oE$ does not vanish identically.\footnote{This assumption can 
be removed at the price of an additional evaluation/interpolation step. 
The details will be covered in the full version of this paper.} In 
which case, it indeed suffices to define $\cH(s,u_0)$ as $\det M_\oE$ itself. 

By construction, there will be exactly $\cM(E)$ rows of 
$M_\oE$ involving coefficients of $f_{m+1}$ and 
exactly $S(\oE)-\cM(E)$ rows involving the parameter $s$. 
Furthermore, the coefficients of $F^*$ can all be assumed to 
be $1$, by Main Theorem 3 of \cite{gcp}.  
So by expanding $\cH(s,u_0)$ in terms of minors, and using 
Hadamard's inequality \cite{mignotte}, it easily follows that 
the coefficient of $t^i$ in $h(t)$ has absolute value at most \begin{tiny}$\left(\begin{mymatrix} \cM(E)\\ i 
\end{mymatrix}\right)$\end{tiny}$\left(2n^2\cM(E)^2\right)^{n\cM(E)}
\left(c\sqrt{k}\right)^{S(\oE)-\cM(E)}$. 
So by Stirling's formula we can conclude the case $m\!=\!n$. \qed 

\begin{rem} 
\label{rem:special} 
The main trick in the above algorithm is to use a 
generic linear form to project the roots of $F$ onto $\C$. This is 
done so as to leave the underlying coordinate rings 
isomorphic over $\Q$. This motivated our choice of $f_{m+1}$ above. 
In particular, the roots of $h$ are exactly the image of 
the roots of $F$ under the chosen linear form. However, this technique 
is quite general and other choices of $f_{m+1}$ are possible and quite 
useful. For instance, picking $f_{m+1}\!=\!u_0+u_ix_i$ 
amounts to projecting the roots onto the $x_i$-axis. 
\end{rem} 

\noindent
{\bf The case $\mathbf{m\!>\!n}$:} 
We first note that the zero set of $F$, now considered as an 
$m\times m$ polynomial system, consists of a set of disjoint 
linear subvarieties. Under the natural embedding of $\C^n 
\hookrightarrow \C^m$ (into the first $n$ coordinates), 
these linear subspaces intersect $\C^n$ precisely in the 
roots of $F$. Furthermore, these linear subspaces never 
intersect within the toric compactification $\cT$ corresponding to 
the polytope $\sum^{m+1}_{i=1} \conv(E_i)$. In particular, 
there is a one to one correspondence between these 
linear subspaces of $\cT$ and the points of $Z$.  

To solve $F$, it thus suffices to find a point 
in every irreducible component of the zero set of $F$ in $\cT$, 
and then project these points onto $\Cn$. Main Theorem 1 of \cite{gcp} 
can do exactly this for us. As for the complexity and 
coefficient analysis, this proceeds almost the same as the 
case $m\!=\!n$ solved above. The only exception is that 
we will need only $n$ of the polynomials $h_1,\ldots,h_m$.  
Our complexity and size bounds then follow immediately from
the $m\!=\!n$ case. \qed  

\section{Coefficient Bounds and Prime Distribution: Proving Main 
Theorem \ref{main:height} and Main Theorem \ref{main:start}} 
\label{sec:proof4} 
 
\noindent
{\bf Proof of Main Theorem \ref{main:height}:}  
The bounds on the coefficients of $h$ derived in 
the proof of Main Theorem \ref{main:solve} (which are 
more precise than those stated in the theorem itself) 
are general enough that we can vary $u_1,\ldots,u_n$ 
to our advantage. In particular, setting all $u_j$ to $0$ 
except for $u_i$, we obtain that any $i^\thth$ coordinate 
of a root of $F$ must be a root of $h$. By 
\cite[theorem 4.2, (viii)]{mignotte}, our root size bound then 
follows immediately. \qed 

\noindent
{\bf Proof of Main Theorem \ref{main:start}:} 
By Main Theorem \ref{main:solve} 
(see also remark \ref{rem:special}), 
we obtain a polynomial $h$ 
with the following properties: (a) $\deg h\!\leq\!\cM(E)$, 
(b) the coefficients of $h$ are integers of   
size polynomial in $S(\oE)$ and the size of $F$, 
and (c) $F$ has a rational root iff $h$ has a rational root. 

Main Theorem \ref{main:start} then follows immediately from the main theorem 
of \cite{weinberger}, which is essentially just a univariate version 
of Main Theorem \ref{main:start}. In particular, our $A$ is the same $A$ as 
that of Weinberger, and the quantity $C_F$ is 
just the logarithm of the discriminant of $h$. By 
Main Theorem \ref{main:height} (and an application of the 
Jacobian), $C_F$ has size polynomial in $S(\oE)$ and the size of $F$. 
\qed 

\section{Acknowledgements} 
The author thanks Felipe Cucker for pointing out the 
excellent references \cite{renegar,bpr}, implicitly containing the 
fact that assertion (1) of Main Theorem 1 could already be done in 
$\mathbf{PSPACE}$ (albeit without monomially sensitive complexity bounds). 

\footnotesize
\bibliographystyle{acm}

\end{document}